\documentstyle[12pt]{article}
\newcommand{\nc}{\newcommand}  \nc{\ov}{\over} \nc{\cd}{\cdots}
\nc{\be}{\begin{equation}} \nc{\ee}{\end{equation}}
\nc{\bZ}{{\bf Z}} \nc{\dl}{\delta} \nc{\la}{\lambda}
\nc{\inv}{^{-1}}  \nc{\ph}{\phi} \nc{\ps}{\psi}
\nc{\twotwo}[4]{\left(\begin{array}{cc}#1&#2\\#3&#4\end{array}\right)}
\nc{\twoone}[2]{\left(\begin{array}{c}#1\\#2\end{array}\right)}
\nc{\cP}{{\cal {P}}} \renewcommand{\sp}{\vskip1ex} \nc{\noi}{\noindent}
\newcommand{\ba}{\begin{array}}
\newcommand{\ea}{\end{array}}

\newcommand{\bt}{\beta}

\renewcommand{\tau}{{\cal D}}

\newcommand{\bea}{\begin{eqnarray}}
\newcommand{\eea}{\end{eqnarray}}

\newcommand{\al}{\alpha}

\newcommand{\Ga}{\Gamma}
\newcommand{\ga}{\gamma}

\newcommand{\pr}{\prime}
\newcommand{\ot}{\frac{1}{2}}
\newcommand{\tp}{\frac{\theta}{\pi}}
\newcommand{\pg}{\frac{i\pi}{g}}
\newcommand{\xg}{\frac{i\xi}{2g}}
\begin{document}\vspace*{3cm}

\begin{center}{\bf A NOTE ON WIENER-HOPF DETERMINANTS\\
    AND THE BORODIN-OKOUNKOV IDENTITY}\end{center}\vskip.1ex
\begin{center}{\bf Estelle L. Basor\footnote{Supported by National
Science Foundation grant DMS--9970879 and EPSRC Grant GR/N35281.} and Yang Chen}
\end{center}\vskip.1ex

\begin{quotation}{\small
The continuous analogue of a Toeplitz determinant
identity for Wiener-Hopf operators is proved. 
An example which arises from random
matrix theory is studied and an error term for the 
asymptotics of the determinant
is computed.}\end{quotation} \vskip.2ex

\section{Wiener-Hopf Determinants}

Recently, a beautiful identity due to Borodin and Okounkov
was proved for Toeplitz determinants 
which shows how one can write a Toeplitz determinant as a 
Fredholm determinant. In this note we generalize this to the 
Wiener-Hopf case. The proof in the Wiener-Hopf case follows 
identically with the second one given in \cite{BW}. We include it 
here for completeness sake and because the nature of the identity is 
slightly different in the continuous verses discrete convolution 
setting. 

In the 
Wiener-Hopf case we begin with a Fredholm determinant on a finite
interval and then show how this can be written as a Fredholm 
determinant of an operator defined on $L^{2}$ of a half-line.
The point is the second operator has a very ``small'' kernel and
thus higher order approximations (as a function of the length of the
finite interval) can be found. 

We now state the analogue of the identity and then apply it to
a particular case to show how error estimates can be computed. In the 
future we hope to refine the estimates given here, apply this identity 
to other important examples, and also extend it to other operators.
 
We consider the Fredholm determinant of the finite Wiener-Hopf 
operator
\be
\det(I-{\cal K}_{[0,\al]}),
\ee
where the ${\cal K}$ acts on $L^2(0,\al)$ and has kernel
 ${\cal K}(x-y) $
with ${\cal K}$ given by the Fourier transform of a function $F$, 
i.e.
\be
F(\xi)=\int_{-\infty}^{\infty}{\cal K}(x)e^{i\xi x}dx.
\ee
 
The continuous analogue of the Borodin-Okounkov
identity is given, under appropriate conditions by the formula
\be
\det(I-{\cal K}_{[0,\al]})=Z\;e^{c\al}\det(I-{\cal L}_{[\al,\infty)}),
\ee
where
\be
c:=\int_{-\infty}^{\infty}\ln(1-F(\xi))\frac{d\xi}{2\pi},
\ee
  ${\cal L}$ is an operator acting on $L^2[\al,\infty)$ 
with kernel,
\be
{\cal L}(x,y)=
\int_{0}^{\infty}\left(\frac{\phi_{-}}{\phi_{+}}
-1\right)_{x+z}
\left(\frac{\phi_{+}}{\phi_{-}}-1\right)_{-z-y}dz,
\ee
and $Z$ is a certain constant whose value will be defined shortly.
Here $\phi_x$ is inverse Fourier transform of 
$\phi(\xi),$ and $\phi_{\pm}(\xi)$ are the Wiener-Hopf
factors satisfying;
\be
1-F(\xi)=\phi_{+}(\xi)\phi_{-}(\xi),\quad \xi\in {\bf R},
\ee
and also satisfying the condition that the functions
$\phi_{\pm}$ when extended away from ${\bf R}$ are
analytic in the upper and lower half-plane respectively.
So that this makes sense and our proof is valid we require that 
$F$ is bounded and in $L^{1}(\bf R)$ and that $\cal K$ is also
in $L^{1}(\bf R)$ and satisfies $\int_{-\infty}^{\infty}|x|{\cal|K|} (x) 
dx < \infty.$ 
To achieve an unique factorization, $1-F$ must have
index zero, be bounded away from zero and we assume that 
$\phi_{\pm}$ are one at $\pm \infty$. These conditions guarantee
that the operator $A - I$ in the proof below 
is trace class and also that all the integrals defined are finite. 
For details see \cite{BS}.

In this section we include the proof of the identity as promised. In 
the next section we apply it to a particular example that arises in 
random matrix theory \cite{Chen}. In this example 
 \bea
{\cal K}(x-y):=\frac{g\sin\pi (x-y)}{\pi\sinh g(x-y)},
\quad x,y\in {\bf R}, \quad g>0.
\eea
We show that as $\alpha \rightarrow \infty,$
\bea
\det(I-{\cal L}_{[\al,\infty)})\sim
1-C(g){\rm e}^{-2g\al(1-\theta/\pi)},
\eea
where $C(g)$ is a completely determined constant and $\cos \theta := 
e^{-\pi^{2}/g},\,\,\,\, 0 < \theta < \pi/2.$ This is a refinement of the 
classical Szeg\"{o}-Kac-Widom Theorem. Previous attempts at the 
refinement computed the exponential term but not the constant. In 
principal, as the reader will see, higher order terms can also be found
using the method outlined in section two.

Here is a proof of the above identity in the Wiener-Hopf case. As 
already stated this proof follows from \cite{BW}. The interested reader should 
also note another slightly different and very elegant proof given in \cite{Bo}.
Denote by $P_{\alpha}$ the orthogonal projection of $L^2(0,\infty)$
onto $L^2(0,\al),$ $Q_{\alpha} = I - P_{\alpha},$ and $P$ the 
orthogonal projection of $L^2(-\infty,\infty)$ onto 
$L^2(0,\infty).$
Also define $W(\phi)$ to be $P M_{\phi} P$ and $W_{\alpha}(\phi)$ to 
be $P_{\alpha} M_{\phi} P_{\alpha}$ where $M_{\phi}$ is multiplication 
by $\phi.$  Note that by using Fourier transforms it can be proved 
that $W_{\alpha}(\phi)$ is unitarily equivalent to $I-{\cal K}_{[0,\al]}$ 
with $\phi = 1 - F.$  It is straight-forward to check that 
 \[
P_{\alpha}\,W(\phi_{+}) = P_{\alpha}\,W(\phi_{+})\,P_{\alpha},\ \ \ 
W(\phi_{-})\,P_{\alpha} = P_{\alpha}\,W(\phi_{-})\,P_{\alpha},\]
and 
\[ W(\phi_{+})W(\phi_{+}^{-1}) = I,\ \ \
W(\phi_{-})W(\phi_{-}^{-1}) = I. \]
Using the above, we can write \footnote{It is an easy general fact that
if 
$\psi_1\in \overline{H^{\infty}}$ or $\psi_2\in H^{\infty}$ then
$W(\psi_1\psi_2)
=W(\psi_1)W(\psi_2)$. In particular $W(\phi_{\pm})$ are invertible with
inverses
$W(\phi_{\pm}^{-1})$. Recall that $H^{\infty}$ consists of all $\psi\in
L^{\infty}$ such that the Fourier transform of $\psi$ vanishes on the 
negative real axis.}
 
\[P_\alpha\,W(\phi)\,P_\alpha
=P_{\alpha}\,W(\phi_{+})\,W(\phi_{+}^{-1})\,W(\phi)\,
W(\phi_{-}^{-1})\,W(\phi_{-})\,P_{\alpha}\]
\[=P_{\alpha}\,W(\phi_{+})\,P_{\alpha}\,W(\phi_{+}^{-1})\,W(\phi)\,
W(\phi_{-}^{-1})\,P_{\alpha}\,W(\phi_{-})\,P_{\alpha}.\]

Now it can be shown that the product of the determinants of  
$P_{\alpha}W(\phi_{\pm})\,P_{\alpha}$ are equal to $e^{c\alpha}$, (see
\cite{BS}, section 10.79).
Thus  to compute $\det(I-{\cal K}_{[0,\al]}) $ we need to consider 
$P_{\alpha}\,W(\phi_{+}^{-1})\,W(\phi)\,W(\phi_{-}^{-1})\,P_{\alpha}$.

Set 
\[W(\phi_{+}^{-1})\,W(\phi)\,
W(\phi_{-}^{-1}) = A.\]
Notice that the determinant of $P_{\alpha}AP_{\alpha}$
equals $\det\,(P_{\alpha}AP_{\alpha}+Q_\alpha)$.   Now $A$ is
invertible 
and differs from $I$ by a trace class operator \cite{BS}. Therefore
\[\det\,( P_{\alpha}A P_{\alpha}+Q_{\alpha})=\det A \,\det\,(A^{-1} P_{\alpha}A
P_{\alpha}+A^{-1}Q_{\alpha})\]
 \[= \det A \,\det\,(A^{-1}( I-Q_{\alpha})A P_{\alpha}+A^{-1}Q_{\alpha}) 
 = \det A \,\det\,(P_{\alpha}-A^{-1}Q_{\alpha}A P_{\alpha}
+A^{-1}Q_{\alpha})\]
 \[= \det A \,\det\,(P_{\alpha} + A^{-1}Q_{\alpha})\,\det\,(I - Q_{\alpha}A P_{\alpha}),\]
 since  $P_{\alpha}Q_\alpha=0$. The determinant of the operator on the right
equals one,
 again since  $P_{\alpha}Q_\alpha=0$. Moreover
 \[\det\,(P_{\alpha}+A^{-1} Q_{\alpha}Q_{\alpha})=\det\, (I-(I-A^{-1})Q_{\alpha})=\det
\,(I-Q_{\alpha}(I-A^{-1})Q_{\alpha}).\]
 We have shown 
 \be \det(I-{\cal K}_{[0,\al]}) =\det A \, \det \,
 (I-Q_{\alpha}(I-A^{-1})Q_{\alpha}).\ee
 It remains to show that this is the same as (1.3). First, $A$
is similar 
 via the invertible operator $W(\phi_+)$ to 
$W(\phi)\,W(\phi_{-}^{-1})\,W(\phi_{+}^{-1}).$ Therefore
\be\det
A=\det\,W(\phi)\,W(\phi_{-}^{-1})\,W(\phi_{+}^{-1})=\det\,W(\phi)\,W(\phi^{-1}).
\label{detA}\ee
This is a representation of the constant $Z$\footnote{It is 
known that the constant $Z$ can also be expressed as 
$\exp\int_{0}^{\infty} z \ln (1 - F)_{z} \ln (1 - F)_{-z}dz$, (see 
\cite{BS} section 10.79).} in 
(1.3).
 Next
 \be
A^{-1}=W(\phi_-)\,W(\phi)^{-1}\,W(\phi_+)=
W(\phi_-)\,W(\phi_+^{-1})\,W(\phi_-^{-1})\,W(\phi_+)
\ee
\be
 =W(\phi_-/\phi_+)\,W(\phi_+/\phi_-). \ee
 Because $\phi_-/\phi_+$ and $\phi_+/\phi_-$ are reciprocals it follows 
 easily from the 
 algebra properties of our operators that the determinant of
 \[
 (I-Q_{\alpha}(I-A^{-1})Q_{\alpha})
 \] is the same as the determinant given in the right-hand side of 
 (3) and this completes the proof.

\section{An Example}

Next we turn to the example cited in the first section.
Consider (1.1)
with
 \bea
{\cal K}(x-y):=\frac{g\sin\pi (x-y)}{\pi\sinh g(x-y)},
\quad x,y\in {\bf R}, \quad g>0.
\eea
The Fourier transform is given by,
\be
F(\xi)=\frac{\sinh(\pi^2/g)}{\cosh(\pi^2/g)+\cosh(\pi\xi/g)}.
\ee
It is clear that $F$ satisfies our hypothesis.
Writing, $1-F(\xi)={\rm exp}(\psi(\xi)),$
and
\be
\psi(\xi)=\int_{0}^{\infty}\psi_t\;e^{i\xi t}dt+
\int_{-\infty}^{0}\psi_t\;e^{i\xi t}dt,
\ee
where 
\be
\psi_t=\int_{-\infty}^{\infty}\ln(1-F(\xi))e^{-i\xi t}
\frac{d\xi}{2\pi}
\ee
yields
\bea
\frac{\phi_{\pm}(\xi)}{\phi_{\mp}(\xi)}
={\rm exp}\left(\pm\,i\Phi(\xi)\right),
\eea
with 
\bea
\Phi(\xi)=2i\int_{0}^{\infty}\psi_t\sin(\xi t)dt.
\eea
A simple calculation shows that
\bea
{\cal L}(x,y)=\int_{0}^{\infty}f(x+z)f(z+y)dz,
\eea
where
\bea
f(x):=\left(\frac{\phi_{+}}{\phi_{-}}-1\right)_{-x}
=\int_{-\infty}^{\infty}
\left(e^{i\Phi(\xi)}-1\right)e^{i\xi x}\frac{d\xi}{2\pi}
=\left(\frac{\phi_{-}}{\phi_{+}}-1\right)_x.
\eea
 To find $f(t)$ we need to compute
$\psi_t$ and $\Phi(\xi).$ With $a:=e^{-\pi^2/g},$ 
a calculation gives,
\bea
\psi_t
=\frac{\cos(\pi t)-\cosh(g\theta t/\pi)}{t\sinh(gt)},
\quad a=:\cos\theta,\;\;0<\theta<\pi/2.
\eea
It is a bit difficult to compute $\Phi(\xi)$ directly, so instead of finding $\Phi(\xi)$
we attempt to find its derivative. 
Using (2.18) and (2.21) and we find, that
\bea
\Phi(\xi) = 2i\int_{0}^{\infty}\left(\frac{\cos(\pi t)-\cosh(g\theta 
t/\pi)}{t\sinh(gt)}\right)\sin(\xi t)dt.
\eea
Next using the above and formula 3.524.5 of \cite{GR}, gives the integral as
a limit as $\mu\to 1$ of the sum of $\zeta(\mu,(1-\beta/\gamma)/2)$ and 
$\zeta(\mu,(1+\beta/\gamma)/2).$ Then using the definition of the $\zeta$ function, 
rearranging the sums, and finally letting $\mu\to 1$ we obtain
\bea
2g\Phi^{\pr}(\xi)
&=&\Psi\left(\ot(1+\tp)-\xg\right)
-\Psi\left(\ot+\frac{i(\pi+\xi)}{2g}\right)
\nonumber\\
&+&\Psi\left(\ot\left(1+\tp\right)+\xg\right)
-\Psi\left(\ot+\frac{i\left(\pi-\xi\right)}{2g}\right)
\nonumber\\
&+&\Psi\left(\ot\left(1-\tp\right)+\xg\right)
-\Psi\left(\ot-\frac{i\left(\pi+\xi\right)}{2g}\right)
\nonumber\\
&+&\Psi\left(\ot\left(1-\tp\right)-\xg\right)
-\Psi\left(\ot-\frac{i\left(\pi-\xi\right)}{2g}\right),
\eea
where $\Psi$ is the di-gamma function. Integrating (2.23)
with respect to $\xi$ with the initial condition
$\Phi(0)=0,$ gives, 
\bea
\mbox{} {\rm exp}(i\Phi(\xi))&=&\frac{\Ga\left(\ot(1+\tp)+\xg\right)}
{\Ga\left(\ot(1+\tp)-\xg\right)}
\frac{\Ga\left(\ot(1-\tp)+\xg\right)}
{\Ga\left(\ot(1-\tp)-\xg\right)}\nonumber\\
&\;&
\times \frac{\Ga\left(\ot(1+\pg)-\xg\right)}
{\Ga(\ot\left(1-\pg)+\xg\right)}
\frac{\Ga\left(\ot(1-\pg)-\xg\right)}
{\Ga(\ot\left(1+\pg)+\xg\right)}.
\eea
We need now to compute (2.20). Note that the integrand vanishes 
when $\xi$ is zero. Furthermore, by using the asymptotics of the Gamma 
functions, it is easily seen that the integral converges conditionally.
Putting $s=\frac{i\xi}{2g}$ and $z=\exp(-2g x)$, in (2.20),
 we find,
\bea
\frac{f(x)}{2g}=\int_{-i\infty}^{i\infty}\left(\frac{\Ga\left(a+s\right)}
{\Ga\left(a-s\right)}
\frac{\Ga\left(1-a+s\right)}
{\Ga\left(1-a-s\right)}
\frac{\Ga\left(b-s\right)}
{\Ga\left(1-b+s\right)}
\frac{\Ga\left(1-b-s\right)}
{\Ga\left(b+s\right)}-1\right)z^{s}\frac{ds}{2\pi i},
\eea
where 
 \bea
\mbox{} a&=&\ot\left(1+\tp\right),\nonumber\\
\mbox{} b&=&\ot\left(1+\pg\right).
\nonumber
\eea
To compute the integral (2.25) we take a sequence of contours 
consisting of the line segment from $-iR_{n}$ to $iR_{n}$ and the 
semi-circle of radius $R_{n}$ in the left-half plane centered at the 
origin. Here $R_{n}=-n+\delta,$ where $\delta$ is any fixed constant 
satisfying $0<\delta<1/4.$ Now computing the 
residues of the Gamma functions in the left-half plane we obtain
after taking the limit $n\to\infty$ the sum of two hypergeometric 
functions $\;_{4}F_{3};$ 
\bea
\mbox{} \frac{f(x)}{2g}&=&
\left(\frac{\Ga(\tp)\Ga(1-\ot(\tp-\pg))\Ga(1-\ot(\tp+\pg))}
{\Ga(1-\tp)\Ga(\ot(\tp+\pg))\Ga(\ot(\tp-\pg))}\right.\\
\mbox{}&\;&\left.\times{\rm e}^{-g\left(1-\tp\right)\al}
\;_4F_3\left(\;^{\bt_1,\,\bt_2,\,\bt_2,\,\bt_1}
_{1-\tp,\,1,\,1-\tp};\exp\left(-2g\al\right)\right)\right)\nonumber\\
\mbox{} &+&
\left(\frac{\Ga(-\tp)\Ga(1+\ot(\tp-\pg))\Ga(1+\ot(\tp+\pg))}
{\Ga(1+\tp)\Ga(\ot(\pg-\tp))\Ga(\ot(-\pg-\tp))}\right.
\nonumber\\
\mbox{}&\;&\left.\times{\rm e}^{-g\left(1+\tp\right)\al}
\;_4F_3\left(\;^{\ga_1,\,\ga_2,\,\ga_2,\,\ga_1}
_{1+\tp,\,1,\,1+\tp};
\exp\left(-2g\al \right)\right)\right),\nonumber
\eea
where
\bea
\mbox{}\bt_1&=&1-\ot\left(\tp-\pg\right),\;\;
\bt_2=1-\ot\left(\pg+\tp\right)
\nonumber\\
\mbox{}\ga_{1}&=&1+\ot\left(\tp-\pg\right),\;\;
\ga_{2}=1+\ot\left(\tp+\pg\right).\nonumber
\eea
As $\al\to\infty,$
\bea
{\rm tr}{\cal L}_{[\al,\infty)}\sim C(g)
\exp\left(-2g\left(1-\tp\right)\al\right),
\eea
where
\bea
C(g):=\left(\frac{\Ga(\tp)\Ga(1-\ot(\tp-\pg))\Ga(1-\ot(\tp+\pg))}{(1-\tp)\Ga(1-\tp)\Ga(\ot(\tp+\pg))
\Ga(\ot(\tp-\pg))}\right)^2.
\eea
We have obtained (2.27) from (2.19) and the trace. Note that
since $0<\theta/\pi<1/2$, the leading exponential
factor, $\exp(-g(1+\theta/\pi)\al)$, in the second term 
of (2.26) will tend to zero faster than the first.
Furthermore, since 
\bea
\;_4F_3
\left(^{A_1,...,A_4}_{B_1,...,B_3};z\right)
=1+\sum_{n=0}^{\infty}\delta_nz^n,\quad |z|<1,\nonumber
\eea
where the above series converges absolutely in the unit circle,
we see that only the exponential factor in the first term of
(2.26) will be relevant in the computation of 
${\rm tr}{\cal L}_{[\al,\infty)}$ for large $\al.$ 
 
Here we give an expression for 
${\rm tr}{\cal L}_{[\al,\infty)}^k,$ where $k$ 
is any positive integer and with this determine 
the higher order terms for large $\al.$ We can compute
$\det(I-{\cal L}_{[\al,\infty)})$ by computing traces.
The justification for this fact follows from our estimate for
$f(x)$ when $x$ is large. It shows that the function is exponentially 
small. Hence when $x$ is larger than $\alpha$ it follows that the 
operator ${\cal L}_{[\al,\infty)}$ has norm smaller than one.

To compute the traces we note that from
the expression of ${\cal L}(x,y)$ given by (2.19) 
and by shifting the interval of integration, $[\al,\infty)$, 
in the trace to $[0,\infty)$,
\bea
{\rm tr}{\cal L}_{[\al,\infty)}^k
=\int_{0}^{\infty}...\int_{0}^{\infty}
f(x_1+x_2+\al)...
f(x_{2k}+x_1+\al)dx_1...dx_{2k}.
\eea
For fixed $x,\,y\,>0$, as $\al\to\infty,$
\bea
f(x+y+\al)\sim 2g
\frac{\Ga(\tp)\Ga(1+\ot(\tp-\pg))\Ga(1-\ot(\tp+\pg))}
{\Ga(1-\tp)\Ga(\ot(\tp+\pg))\Ga(\ot(\tp-\pg))}
{\rm e}^{-g(1-\tp)(x+y+\al)},
\eea
where the second exponentially decaying term in (2.26) has 
been discarded.
Integrating (2.29) with (2.30) we find, as $\al\to\infty,$
\bea
{\rm tr}{\cal L}_{[\al,\infty)}^k
\sim C(g)^k\exp\left(-2g\left(1-\tp\right)k\al\right).
\eea
Therefore,
\bea
-\ln\det(I-{\cal L}_{[\al,\infty)})
&=&\sum_{k=1}^{\infty}\frac{{\rm tr}{\cal L}^k_{[\al,\infty)}}{k},
\nonumber\\
&\sim&\sum_{k=1}^{\infty}\frac{C^k(g)}{k}{\rm e}^{-2g(1-\tp)\al k}.\nonumber
\eea
Finally,
\bea
\det(I-{\cal L}_{[\al,\infty)})\sim
1-C(g){\rm e}^{-2g\al(1-\theta/\pi)},\quad \al\to\infty.
\eea
In random matrix theory $E(\al):=\det(I-{\cal K}_{[0,\al]})$ 
is the probability that an interval $[-\al/2,\al/2]$ 
(after a suitable scaling) is free of eigenvalues. Putting 
$g=0,$ in (2.1), we have the {\it sine} kernel,
$\frac{\sin(\pi(x-y))}{\pi(x-y)}.$ It was shown in 
\cite{JMMS} that $R(\al):=-\frac{d\ln E(\al)}{d\al},$ 
satisfies a particular Painl\'eve V equation. For a simpler 
derivation of this and a review of random matrix theory 
see \cite{TW}. Recently, it was shown in \cite{Ni} using the 
theory of \cite{TW1} that $R(\al)$ for $g>0,$ satisfies 
a particular Painl\'eve VI. In a heuristic perturbative
calculation on the Painl\'eve VI both
the exponential decaying terms were was found \cite{Chen1}, 
however, the prefactor, $C(g),$ cannot be determined as
it involves an indefinite integral. It appears that
$C(g)$ can only be obtained from determinant identity, (1.3). 
It can also be shown that the computation of (1.1) with the
kernel given by (2.13) can be reformulated as an equivalent
$2\times 2$ matrix Riemann-Hilbert problem, see \cite{K} 
for a description of this technique. Although as $\al\to\infty$, 
such a technique could reproduce the classical results of 
Akhiezer, Hirschman and Szeg\"o, namely, the $Z$ factor, it is not at 
all clear whether the  higher correction terms given by (2.34) can be found. 
Final remark: the determinant identity, (1.3), does not hold
for $g=0$. Since in this situation, the Fourier transform of
the {\it sine} kernel is the characteristic function of 
$[-\pi/2,\pi/2]$ and the Wiener-Hopf factorization fails.

\hspace{-2em}\begin{tabular}{lll} Department of Mathematics&&Department
of Mathematics
\\ California Polytechnic State University&&Imperial College\\
San Luis Obispo, CA 93407 USA&&180 Queen's gate, London, SW7 2BZ, UK\\
ebasor@calpoly.edu&&y.chen@ic.ac.uk\end{tabular}\sp
\vspace{3ex}

\end{document}